\documentclass[review]{siamart171218}

\usepackage{amsfonts}
\usepackage{graphicx}
\usepackage{algorithm,algpseudocode}
\usepackage{epsfig,amssymb,amsmath,listings,lineno,xcolor,bm,url}
\usepackage{tabulary}
\usepackage{booktabs}
\usepackage[normalem]{ulem}
\usepackage{subfig}
\usepackage{float}
\ifpdf%
  \DeclareGraphicsExtensions{.pdf,.pdf,.png,.jpg}
\else
  \DeclareGraphicsExtensions{.pdf}
\fi

\newtheorem{remark}{Remark}[section]

\newcommand{\TheTitle}{%
  Mixed precision matrix interpolative decompositions for model reduction
}

\newcommand{\TheShortTitle}{%
  Mixed Precision Interpolative Decomposition
}

\newcommand{\TheName}{%
  Alec Dunton, Alyson Fox
}

\newcommand{\TheAddress}{%
  University of Colorado Boulder,
  (\email{alec.dunton@colorado.edu}).
}

\newcommand{\TheFunding}{%
  This work was funded in part by the United States Department of Energy's National Nuclear Security Administration under the Predictive Science Academic Alliance Program II at Stanford University, Grant DE-NA-0002373. This project was also funded by and conducted as a part of the Computing Scholar Program at Lawrence Livermore National Laboratory. This work was performed under the auspices of the U.S. Department of Energy by Lawrence Livermore National  Laboratory under Contract DE-AC52-07NA27344, LLNL-JRNL-813055-DRAFT
}

\nolinenumbers

\author{\TheName\thanks{\TheAddress}}
\title{{\TheTitle}\thanks{\TheFunding}}
\headers{\TheShortTitle}{\TheName}
\ifpdf%
\hypersetup{%
  pdftitle={\TheTitle},
  pdfauthor={\TheName}
}
\fi

\begin{document}
\maketitle


\begin{abstract} Renewed interest in mixed-precision algorithms has emerged due to growing data capacity and bandwidth concerns, as well as the advancement of GPU’s, which enable significant speedup for low precision arithmetic.
In light of this, we propose a mixed-precision algorithm to generate a double-precision accurate matrix interpolative decomposition (ID) approximation under a given set of criteria. Though low precision arithmetic suffers from quicker accumulation of round-off error, for many data-rich applications we nevertheless attain viable approximation accuracies, as the error incurred using low precision arithmetic is dominated by the error inherent to low-rank approximation. We then conduct several simulated numerical tests to demonstrate the efficacy of the algorithms. Finally, we present the application of our algorithms to a problem in model reduction for particle-laden turbulent flow.
\end{abstract}


\section{Background}\label{sec:intro}
Mixed precision algorithms have recently gained popularity for several reasons. In addition to the reduced memory footprint of lower precision formats, new hardware capabilities like GPUs, e.g., NVIDIA TensorCores \cite{markidis2018nvidia}, enable computers to complete arithmetic in low precision at a much faster rate. In particular, dense linear systems may be solved 2 to 4 times faster in half precision than in double precision arithmetic~\cite{haidar2018design,haidar2018harnessing}, while dense matrix-matrix multiplication may be completed up to 10 times faster~\cite{abdelfattah2019fast}. The data movement bottleneck problem, both on and off node, has increasingly become the dominating barrier to exascale computing; pairing the benefits of the new hardware capabilities with the benefit of the reduced data capacity and bandwidth ensures significant speed ups. However, as lower precisions are introduced into algorithms, concerns about floating-point round-off error must be addressed. Despite these concerns, lower precision arithmetic has been incorporated into many numerical linear algebra algorithms, among them iterative refinement~\cite{carson2018accelerating} and QR factorization~\cite{yang2019mixed}.  
We propose that mixed precision can be used to compute low-rank pivoted QR factorizations; the errors accrued in lower precision will be dominated by low-rank approximation error; these decompositions can then be used to form the matrix interpolative decomposition (ID).

IEEE-754 floating point number systems constitute subsets of the reals $F \subset \mathbb{R}$, representing numbers in terms of a base $b \in \mathbb{N}$, precision $t$, significand $m \in \mathbb{N}$ with $0 \leq m \leq b^{t}-1$, and an exponent range $e_{min} \leq e \leq e_{max}$. For every element $y \in F$ we have 
\begin{equation}
    y = \pm m \times b^{e - t}
\end{equation}
The bits used to encode $e$ are referred to as exponent bits, the bits used to represent $m$ the mantissa bits, with one bit reserved for encoding the sign of the element $y \in F$.

In this work, mixed precision arithmetic constitutes combined use of the three IEEE-754 formats: 16 bit binary (half precision), 32 bit binary (single precision), and 64 bit binary (double precision). 
The number of exponent bits, mantissa bits, and unit round-off errors for these three formats are provided in Table~\ref{tab:bits}.

By definition, any two distinct numbers represented in floating point format have some finite spacing between them; rounding will necessarily introduce round-off error to any binary arithmetic operation. We define $fl(\cdot)$ to be the rounding operation, which maps a real number to its closet IEEE-754 representation, and $u = \frac{1}{2} b^{1-t}$ to be the unit round-off error. In order to model round-off error, we follow~\cite{higham2002accuracy}:
\begin{equation}
    fl (x \hspace{3pt} \mbox{\textbf{op}} \hspace{3pt} y) = (x \hspace{3pt} \mbox{\textbf{op}} \hspace{3pt} y)(1 + \delta) , \hspace{3pt}  \vert \delta \vert \leq u ,
\end{equation}
where \textbf{op} represents addition, subtraction, multiplication, or division. 
We now extend this model to successive FLOPs with the following definition. 
\begin{definition}
(Lemma 3.1 in~\cite{higham2002accuracy}.) Let $\vert \delta_i \vert < u$ and $\rho_i \in \{ -1, +1 \}$ for $i = 1, \dots, k$ and $ku < 1$. Then,
\begin{equation}
\prod_{i=1}^k (1 + \delta_i)^{\rho_i} = 1 + \theta^{(k)}, \hspace{6pt} where \hspace{6pt} \vert \theta^{(k)} \vert \leq \frac{ku}{1 - u} =: \gamma_k .
\end{equation}
\end{definition}
The quantity $\theta^{(k)}$ represents accumulation of roundoff error in $k$ successive FLOPs, and is bounded above by $\gamma_k$ to provide notational convenience.
All notations used in this work are provided in Table~\ref{tab:TableOfNotationForMyResearch}.

\begin{table}
\centering,
\begin{tabular}{ |p{3cm}|p{2.7cm}|p{2.7cm}|p{2.7cm}| }
\hline
IEEE Standard & Exponent Bits & Mantissa Bits  & Round-off Error \\
\hline
Half-precision  & 5  & 10 & $4.9 \times 10^{-4}$\\
Single-precision & 8   & 23 &$6.0 \times 10^{-8}$\\
Double-precision & 11  & 52  & $1.1 \times 10^{-16}$\\
\hline
\end{tabular}
\caption{Characteristics of three IEEE-754 float formats}
\label{tab:bits}
\end{table}
\label{sec:notation}
\begin{table}[htbp]
\centering 
\begin{tabular}{r c p{10cm} }
\toprule
$\bm{A}$ & $\triangleq$ & data matrix\\
$\mathcal{I}$ & $\triangleq$ & column index vector\\
$\bm{A}(:,\mathcal{I})$ & $\triangleq$ & column skeleton of data matrix\\ 
$fl(\cdot)$ & $\triangleq$ & floating-point rounding operator \\
$u$ & $\triangleq$ & unit round-off error \\
$u_L$ & $\triangleq$ & low precision round-off error \\
$\gamma_k$ & $\triangleq$ & round-off error accumulated from $k$ successive FLOPs \\
$\bm{A}_D$ & $\triangleq$ & double-precision original data matrix \\
$\sigma_j(\bm{A}_D)$ & $\triangleq$ & $j^{th}$ largest singular value of $\bm{A}_D$\\
$\bm{A}_{D(i,j)}$ & $\triangleq$ & entry of $\bm{A}_D$ in $i^{th}$ row and $j^{th}$ column \\
$\bm{A}_L$ & $\triangleq$ & low-precision data matrix \\
$\hat{\bm{A}}_D$ & $\triangleq$ & double-precision ID approximation\\
$\hat{\bm{A}}_M$ & $\triangleq$ & mixed-precision ID approximation\\
$\hat{\bm{A}}_L$ & $\triangleq$ & low-precision ID approximation\\
$+$ & $\triangleq$ & Moore-Penrose pseudo-inverse \\
\bottomrule
\end{tabular}
\caption{Notation Table}
\label{tab:TableOfNotationForMyResearch}
\end{table}

\subsection{Column ID}
Broadly, low-rank matrix approximations seek to identify {\it factor matrices} $\bm{B} \in \mathbb{R}^{m \times k}$ and $\bm{C} \in \mathbb{R}^{k \times n}$ with $k \ll m,n$ such that $\Vert \bm{A} - \bm{BC} \Vert \ll \epsilon$ for some $0 < \epsilon \ll 1$. This work focuses on one such approximation: the column interpolative decomposition (column ID). We approximate  a matrix $\mathbf{A} \in \mathbb{R}^{m \times n}$ via the column ID as a product of $k$ of its columns $\mathbf{A(:,\mathcal{I}}) \in \mathbb{R}^{m \times k} $, indexed by $\mathcal{I}\subseteq\{1,\dots,n\}$ with $\vert\mathcal{I}\vert=k$, referred to as a {\it column skeleton}, and a {\it coefficient matrix} encoding an interpolation rule on its columns $\bm{a}_i$, which we denote $\mathbf{P} \in \mathbb{R}^{k \times n}$~\cite{cheng2005compression}. In this case, our factor matrix $\bm{B}$ is the column skeleton $\mathbf{A(:,\mathcal{I}})$ and $\bm{C}$ is the coefficient matrix $\mathbf{P}$. Mathematically, the $k$-rank interpolative decomposition is then defined as
\begin{equation}
     \begin{bmatrix}
     \vert & \vert & \vert & \vert \\
     \bm{a}_1 & \bm{a}_2 & \dots & \bm{a}_n \\
     \vert & \vert & \vert & \vert
     \end{bmatrix} 
    \approx \begin{bmatrix}
     \vert & \vert & \vert & \vert \\
     \bm{a}_{i_1} & \bm{a}_{i_2} & \dots & \bm{a}_{i_k}  \\
     \vert & \vert & \vert & \vert
     \end{bmatrix} 
     \bm{P}
     \hspace{4pt}
     ;
\end{equation}
where,
\begin{equation}
\hspace{4pt}
     \bm{a}_i = \sum_{j=1}^k \bm{a}_{i_j}\bm{P}_{ij} 
     \hspace{4pt}
     ;
     \hspace{4pt}
     i_j \in \mathcal{I}
     \hspace{4pt}.
\end{equation}

In order to construct an ID, we first take the input matrix $\bm{A}$ and compute a column pivoted QR decomposition using a modified gram-schmidt scheme~\cite{golub2012matrix,cheng2005compression} so that
\begin{equation}
\bm{AZ} = \bm{QR}, 
\end{equation}
where $\bm{Z}$ is a permutation matrix. To obtain a $k$-rank approximation of $\bm{A}$, we rewrite the QR decomposition in block format;
\begin{equation}
    \bm{Q}\bm{R} = \begin{bmatrix}
    \bm{Q}_{11} & \bm{Q}_{12} \end{bmatrix} 
    \begin{bmatrix}
    \bm{R}_{11} & \bm{R}_{12} \\
    \bm{0} & \bm{R}_{22} 
\end{bmatrix},
\end{equation}
where $\bm{Q}_{11}$ is of dimension $m \times k$, $\bm{R}_{11}$ is $k \times k$, and $\bm{R}_{12}$ is $k \times (n-k)$. We then rearrange the sub-matrices;
\begin{equation}
\label{eqn:idconstruction}
\bm{A} \approx \bm{Q}_{11}\bm{R}_{11} \left[\bm{I} \hspace{3pt} \bm{R}_{11}^{+}\bm{R}_{12} \right]\bm{Z}^{T} 
= \bm{A}(:,\mathcal{I})\left[\bm{I} \hspace{3pt} \bm{R}_{11}^{+}\bm{R}_{12} \right]\bm{Z}^{T}
= \bm{A}(:,\mathcal{I})\bm{P}.
\end{equation}
Hence, $\bm{P} = \left[\bm{I} \hspace{3pt} \bm{R}_{11}^{+}\bm{R}_{12} \right]\bm{Z}^{T}$, where $+$ denotes the Moore-Penrose pseudo-inverse. 
\begin{remark}
\label{rem:Tmatrix}
In practice $\bm{R}_{11}$ is frequently ill-conditioned~\cite{cheng2005compression}, which can lead to accumulation of significant errors. Due to this ill-conditioning/rank-deficiency in $\bm{R}_{11}$, we use the pseudo-inverse to construct the coefficient matrix in (\ref{eqn:idconstruction}).
\end{remark}

\section{Mixed Precision ID}
\begin{algorithm}[t]
\caption{Double Precision Column ID $\bm{A} \approx \bm{A}(:,\mathcal{I})\bm{P}$} \label{alg:doubleid} 
\begin{algorithmic}[1]
\Procedure{MPID}{$\bm{A}$ $\in \mathbb{R}^{m \times n}$}
\State $k \gets$ target rank
\State $\bm{Q}$, $\bm{R}$, $\mathcal{I} \gets$ $MGSQR(\bm{A},k)$ 
\State $\bm{T} \gets (\bm{R}(1:k,1:k))^{+}\bm{R}(1:k,(k+1):n)$
\State $\bm{P}(:,\mathcal{I}) \gets \left[\bm{I}_k \hspace{12pt}\bm{T} \right]$
\State $\mathcal{I} \gets \mathcal{I}(1:k)$
\State $\bm{return}$ $\mathcal{I}, \bm{P}$ set column skeleton to $\bm{A}_D(:,\mathcal{I})$.
\EndProcedure
\end{algorithmic}
\end{algorithm}

\begin{algorithm}[t]
\caption{Mixed/Low Precision Column ID $\bm{A} \approx \bm{A}(:,\mathcal{I})\bm{P}$} \label{alg:mixedid} 
\begin{algorithmic}[1]
\Procedure{MPID}{$\bm{A}$ $\in \mathbb{R}^{m \times n}$}
\State $\bm{A}_{L} \gets$ reduced precision $fl(\bm{A})$
\State $k \gets$ target rank
\State $\bm{Q}$, $\bm{R}$, $\mathcal{I} \gets$ (Mixed/Low) $ MGSQR_{L}(\bm{A}_{L},k)$
\State $\bm{T} \gets (\bm{R}(1:k,1:k))^{+}\bm{R}(1:k,(k+1):n)$
\State $\bm{P}(:,\mathcal{I}) \gets \left[\bm{I}_k \hspace{12pt}\bm{T} \right]$
\State $\mathcal{I} \gets \mathcal{I}(1:k)$
\State $\bm{return}$ $\mathcal{I}, \bm{P}$ set column skeleton to $\bm{A}_D(:,\mathcal{I})$ (Mixed) or $\bm{A}_L(:,\mathcal{I})$ (Low).
\EndProcedure
\end{algorithmic}
\end{algorithm}

\label{sec:main}
\label{sec:MPID}
In order to obtain a mixed precision ID algorithm, we augment the column ID algorithm (Algorithm~\ref{alg:doubleid}) so that our double precision input matrix $\bm{A}_D$ is cast to single or half precision, which we denote $\bm{A}_L$ (presented in Algorithm~\ref{alg:mixedid}). Let $u_L$ be the low precision round-off error with respect to the real numbers. By definition, $\bm{A}_D = \bm{A}_L +\bm{E}$, where $|\bm{E}_{(i,j)}|\leq |\bm{A}_{L(i,j)}|u_L$. Then, the index set $\mathcal{I}$ and coefficient matrix $\bm{P}$ are computed using our modified Gram-Schmidt procedure in low precision arithmetic ($MGSQR_L(\cdot)$ from Line 4 in Algorithm~\ref{alg:mixedid}). The indices and coefficient matrix generated in low precision are then used to construct an ID for the original double precision matrix. $MGSQR(\cdot)$ dominates the computational cost in column ID; low precision arithmetic in this step enables overall algorithm speedup if the implementation of the algorithm is either memory or bandwidth bound. 

One of the outputs of $MGSQR(\cdot)$ is the index set $\mathcal{I}$, corresponding to the selected columns $\bm{A}(:,\mathcal{I})$. If the column selection step is affected in low precision, and the columns chosen by the mixed precision algorithm constitute a worse basis than those chosen in the double precision case, we may see degradation of the final approximation. However, the data we are concerned with representing is riddled with numerical approximation error. Therefore, our methods do not tell us what an optimal basis for explaining the underlying physical phenomenon in the first place. We assert that the error incurred in this step of the algorithm will not be the dominating error under suitable conditions.

As well as the column selection step in $MGSQR(\cdot)$, another source of error is the computation of the coefficient matrix $\bm{P}$, which is computed from the $\bm{R}$ matrix and pivot indices output by $MGSQR(\cdot)$ via a least-squares problem. In particular, the computation of $\bm{P}$ typically entails solving an ill-conditioned linear system (due specifically to the conditioning of $\bm{R}_{11}$ in Equation~\ref{eqn:idconstruction}), and therefore may be susceptible to inaccurate solutions or accumulation of error (see Remark~\ref{rem:Tmatrix}). 

Despite concerns with round-off in lower precision arithmetic, in many relevant applications there is inherent error in a given simulation or set of measurements, e.g., structural uncertainty (due to simplifying assumptions in the model used for simulation), or experimental uncertainty (due to noisy measurements). Consequently, the round-off error incurred in low precision arithmetic may be much smaller than these other sources of error, and therefore acceptable. Moreover, if the speedup in mixed precision is significant, then many application areas would welcome this additional error, provided we can assure that the error due to our methods is small compared to that of the dominating error.


We begin with a lemma outlining properties of the ID.

\begin{lemma}
\label{lemma:id}
(Lemma 1 from~\cite{liberty2007randomized}) There exists a rank $k$ column ID $\bm{A} = \bm{A}(:,\mathcal{I})\bm{P}$ such that, in exact arithmetic,
\begin{enumerate}
    \item $\Vert \bm{P} \Vert_{2} \leq \sqrt{1 + k(n-k)}$,
    \item $\Vert \bm{A} - \bm{A}(:,\mathcal{I})\bm{P} \Vert_2 \leq \sqrt{1 + k(n-k)} \sigma_{k+1}(\bm{A}_D)$,
\end{enumerate}
where $n$ is the column dimension of $\bm{A}$ and $\sigma_{k+1}(\bm{A}_D)$ is the $(k+1)^{th}$ largest singular value of the double precision data matrix $\bm{A}_D$.
\end{lemma}

 As in Table~\ref{tab:TableOfNotationForMyResearch}, we denote our original, ground truth double precision matrix $\bm{A}_D$. The double precision ID approximation to our original matrix generated in double precision is defined to be $\hat{\bm{A}}_D = \bm{A}_D(:,\mathcal{I}_D)\bm{P}_D$; the column skeleton and coefficient matrix are both computed entirely in double precision. The mixed precision ID approximation to our double precision matrix is defined to be $\hat{\bm{A}}_M = \bm{A}_D(:,\mathcal{I}_L)\bm{P}_L$, where $\bm{P}_L$ is the coefficient matrix generated in low precision but we use the double precision columns of $\bm{A}_D$ corresponding to the indices identified by the low precision QR. Finally, the low precision ID approximation to our original matrix is defined as $\hat{\bm{A}}_L = \bm{A}_L(:,\mathcal{I}_L)\bm{P}_L$. In the low precision ID approximation we use both the low precision coefficient matrix (as in mixed precision) and the low precision column skeleton to form our final ID approximation to $\bm{A}_D$.

\begin{remark}
In this work, we are concerned with obtaining approximations to data matrices which are originally stored in double precision, and therefore assume the products $\hat{\bm{A}}_D = \bm{A}_D(:,\mathcal{I}_D)\bm{P}_D$, $\hat{\bm{A}}_M = \bm{A}_D(:,\mathcal{I}_L)\bm{P}_L$, and $\hat{\bm{A}}_L = \bm{A}_L(:,\mathcal{I}_L)\bm{P}_L$ are evaluated in double precision. 
\end{remark}

\section{Numerical Results}\label{sec:num}
\label{sec:numericalresults}
We present numerical results for three datasets which we have artificially generated to exhibit specific singular value profiles. Our datasets feature (1) linearly decaying singular values ($\sigma_i = i^{-1}$ - denoted `Slow' in plots), (2) quadratically decaying singular values ($\sigma_i = i^{-2}$ - denoted `Medium' in plots), and (3) quartically decaying singular values ($\sigma_i = i^{-4}$ - denoted `Fast' in plots). In the test cases where we vary the column dimension of the matrix, we fix the target rank to be $20$, i.e., $k =20$, and the row dimension to be $m= 1000$. Increasing the column dimension amounts to sampling more columns from each of the datasets in Table~\ref{tab:datasets}. For example, if we set the column dimension to $n=300$, this is equivalent to defining our input matrices to be the first $300$ columns of the original data matrices. We enumerate properties of our three test datasets in Table~\ref{tab:datasets}. We provide (1) the dimensions, (2) the ratio of the $50^{th}$ largest singular value of the matrix to the largest singular value, (3) the ratio of the smallest singular value to the largest singular value, and (4) the value range of each dataset. We define the value range to be the ratio of the maximum absolute value entry of a matrix to its minimum absolute value entry. 

We select our low precisions to be 32 bit binary (IEEE-754 single precision) and 16 bit binary (IEEE-754 half precision). We implement our methods in Julia v1.1, in which we cast our variables to either half or single precision, but accumulate error exclusively in single precision. We refer to this as simulated half-precision; the interested reader is referred to Algorithm 1 in~\cite{yang2019mixed} for more details. Accumulating error and conducting arithmetic in single precision may lead to overly optimistic results, though this is precisely the manner in which operations are executed on, e.g., NVIDIA TensorCores~\cite{markidis2018nvidia}. Casting variables to low precision requires significant computation time, so we do not measure runtime in these experiments. We report relative accuracy in terms of the matrix 2-norm. In reporting the accuracy of our schemes, we test two important metrics of our algorithms' performances: how similar of an approximation they generate to that of their double precision counterparts, as well as their approximation accuracy relative to the ground truth double precision matrix.

\subsection{Comparison to double precision approximation}

\begin{table}[h]
\centering
\begin{tabular}{ |p{2cm}|p{2cm}|p{2cm}|p{2cm}|p{2cm}| }
\hline
Dataset &  Dimensions & $\sigma_{50}/\sigma_{1}$ & $\sigma_{n}/\sigma_{1}$ & value range\\
\hline
Slow & $1000 \times 1000$ & $2.0 \times 10^{-2}$ & $1\times 10^{-3}$ & $5.0 \times 10^{6}$\\
Medium & $1000 \times 1000$ & $4.0 \times 10^{-4}$ & $1\times 10^{-6}$ & $1.7 \times 10^{6}$\\
Fast  & $1000 \times 1000$ & $1.6 \times 10^{-7}$   & $1 \times 10^{-12}$ & $4.2 \times 10^{6}$ \\
\hline
\end{tabular}
\caption{Characteristics of three datasets to be used in the numerical experiments in this paper.}
\label{tab:datasets}
\end{table}

\begin{figure}
\centering
\includegraphics[width=0.49\linewidth]{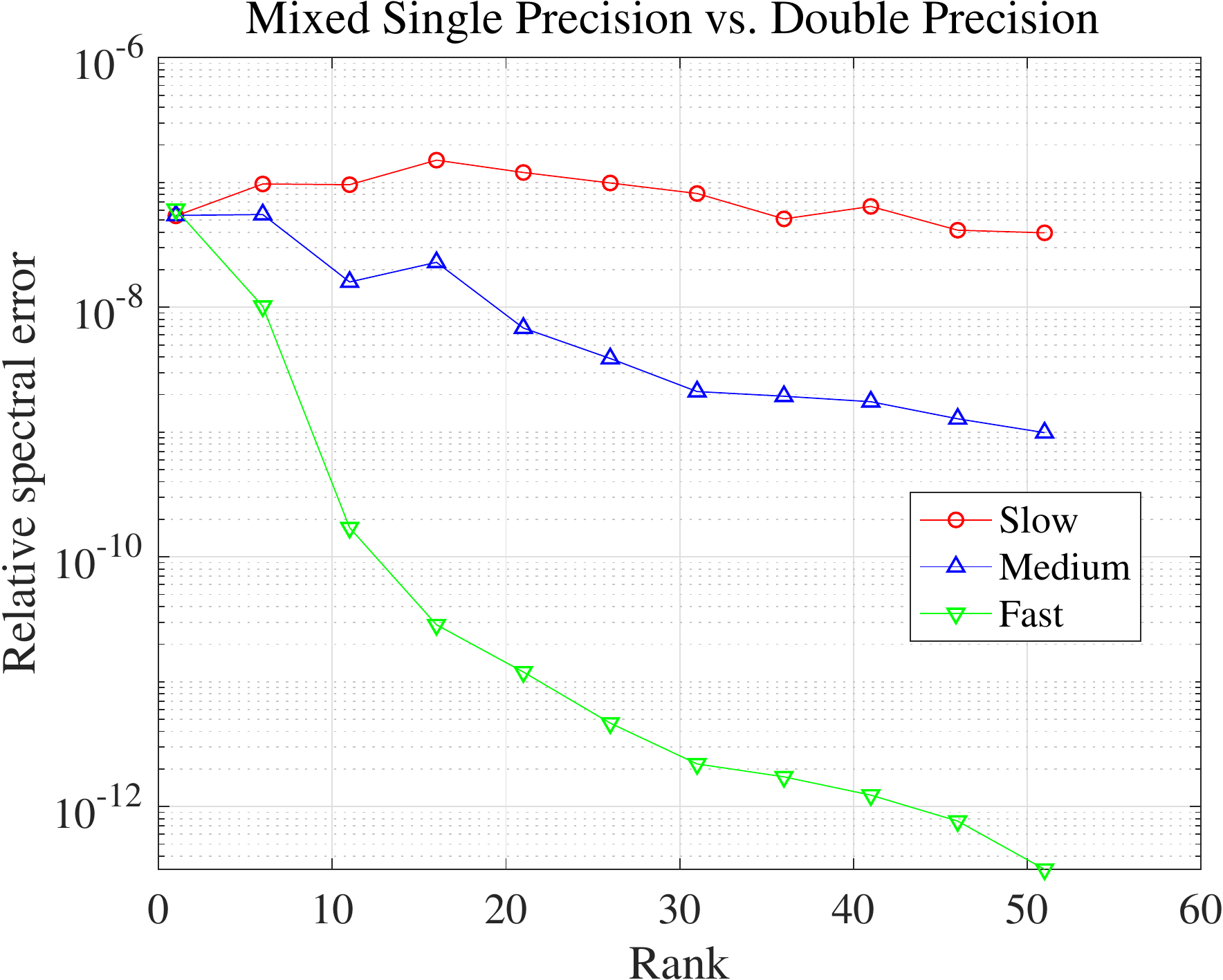}
\includegraphics[width=0.49\linewidth]{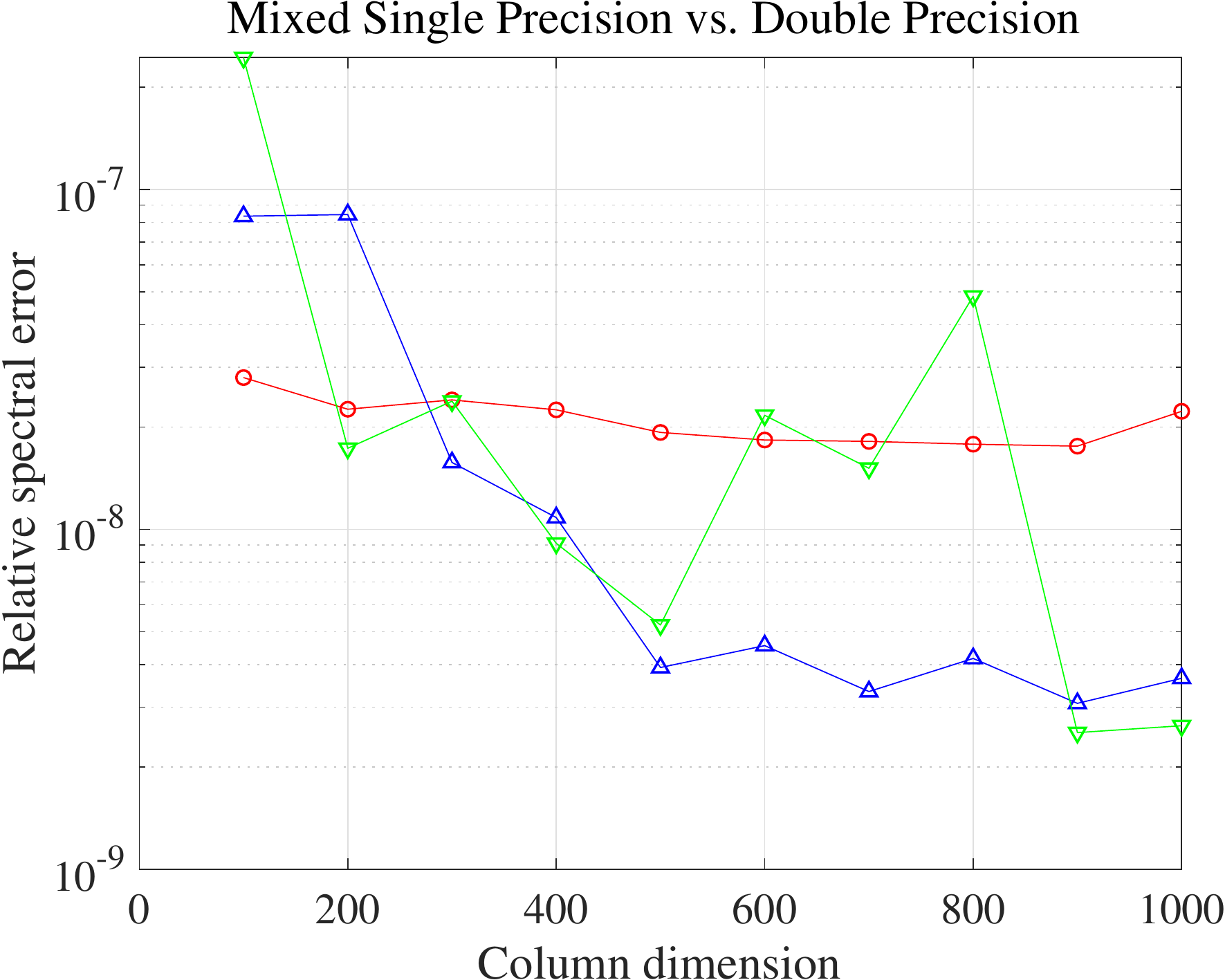}
\caption{Relative spectral error of mixed single precision ID with respect to target rank (left) and column dimension (right).}
\label{fig:mp_dp}
\end{figure}

The value range of each matrix is bounded above by $5 \times 10^6$ for all datasets. Therefore, all entries across all datasets can be represented exactly in single precision, in which unit round-off is $6 \times 10^{-8}$. In light of this, we expect accurate results using our mixed precision ID algorithm with our low precision set to single precision (referred to mixed single precision ID for brevity). We now compare mixed single precision ID to double precision ID. We generate the index set $\mathcal{I}$ in both double and single precision, meaning the indices of the columns chosen in the two precisions could differ. In the left panel of Figure~\ref{fig:mp_dp}, we observe that for the fast and medium decay datasets almost all spectral errors are smaller than single precision unit round-off ($6 \times 10^{-8}$). This indicates that our mixed single precision algorithm performs as well as can be expected in single precision when compared to a double precision algorithm in these two test cases. In the slow decay test case, we also obtain good accuracy, but our errors are larger than unit round-off. 

Relative to the error incurred by low-rank approximation, which is bounded above by $\sqrt{1 + k(n-k)}\sigma_{k+1}(\bm{A}_D)$ (Lemma~\ref{lemma:id}) - roughly two orders of magnitude larger than the singular values of each dataset (see Table~\ref{tab:datasets}) - all measured errors are negligible for the three datasets. More specifically, for the fast decay test case with the largest target rank tested ($k=51$), our low-rank approximation error is bounded above by $3.3 \times 10^{-5}$. This is the best low-rank approximation error across all target rank values and datasets. Therefore, the errors incurred using single precision instead of double precision to construct an ID are all negligible relative to the inherent error in low-rank approximation in all cases.

Increasing the column dimension of our input, we should expect to observe that the relative error of our mixed precision approximation is not impacted. By adding columns to our test matrix, ID is given a larger set from which to construct a rank $20$ basis to approximate the data, although the number of columns being approximated also increases. Therefore, the error accrued by increasing the size of the matrix should be offset by the error due to increased round-off. This is confirmed in the right panel of Figure~\ref{fig:mp_dp}, as we observe that our mixed single precision algorithm is robust with respect to increases in column dimension. For all three datasets, mixed single precision ID achieves error around $10^{-8}$.

\begin{figure}
\centering
\includegraphics[width=0.49\linewidth]{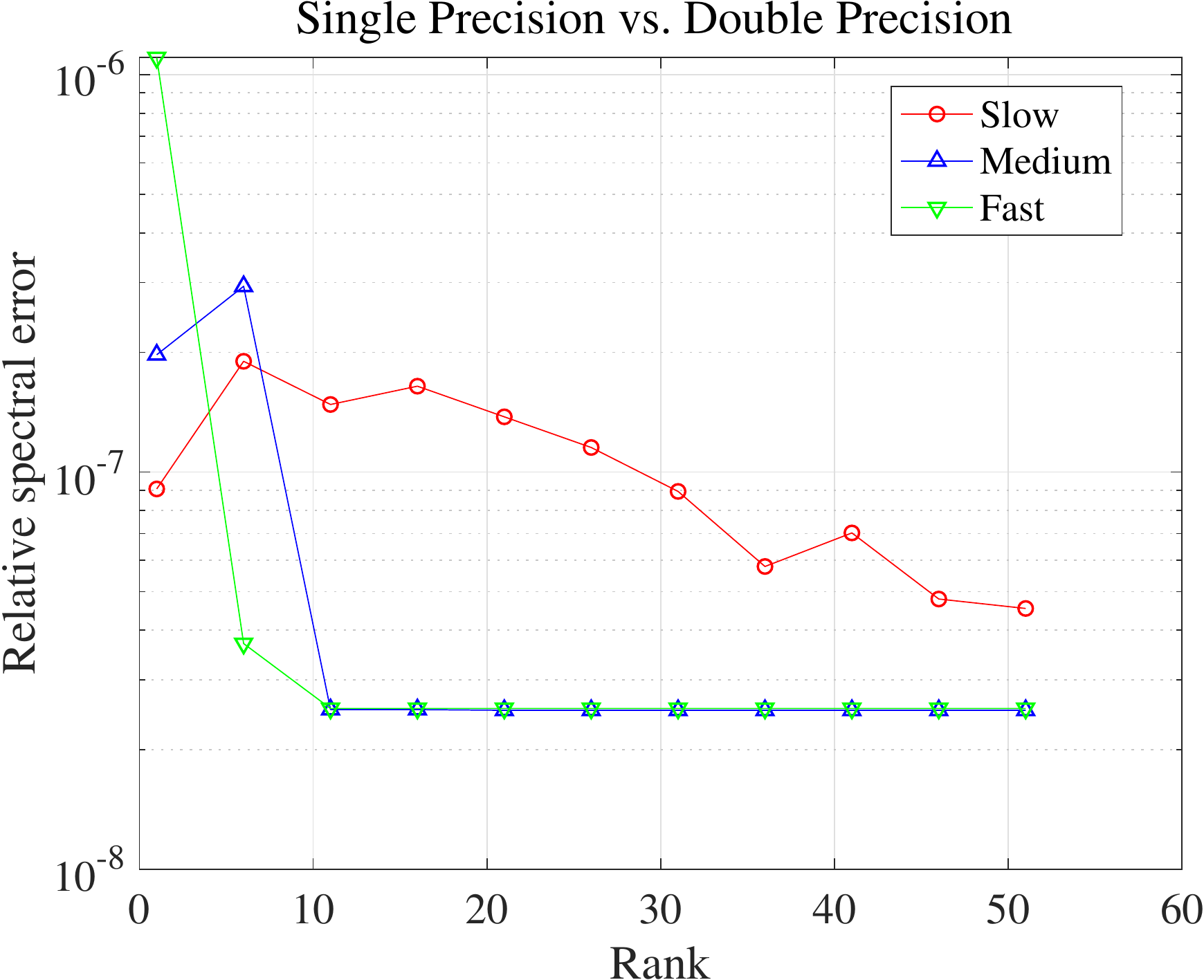}
\includegraphics[width=0.49\linewidth]{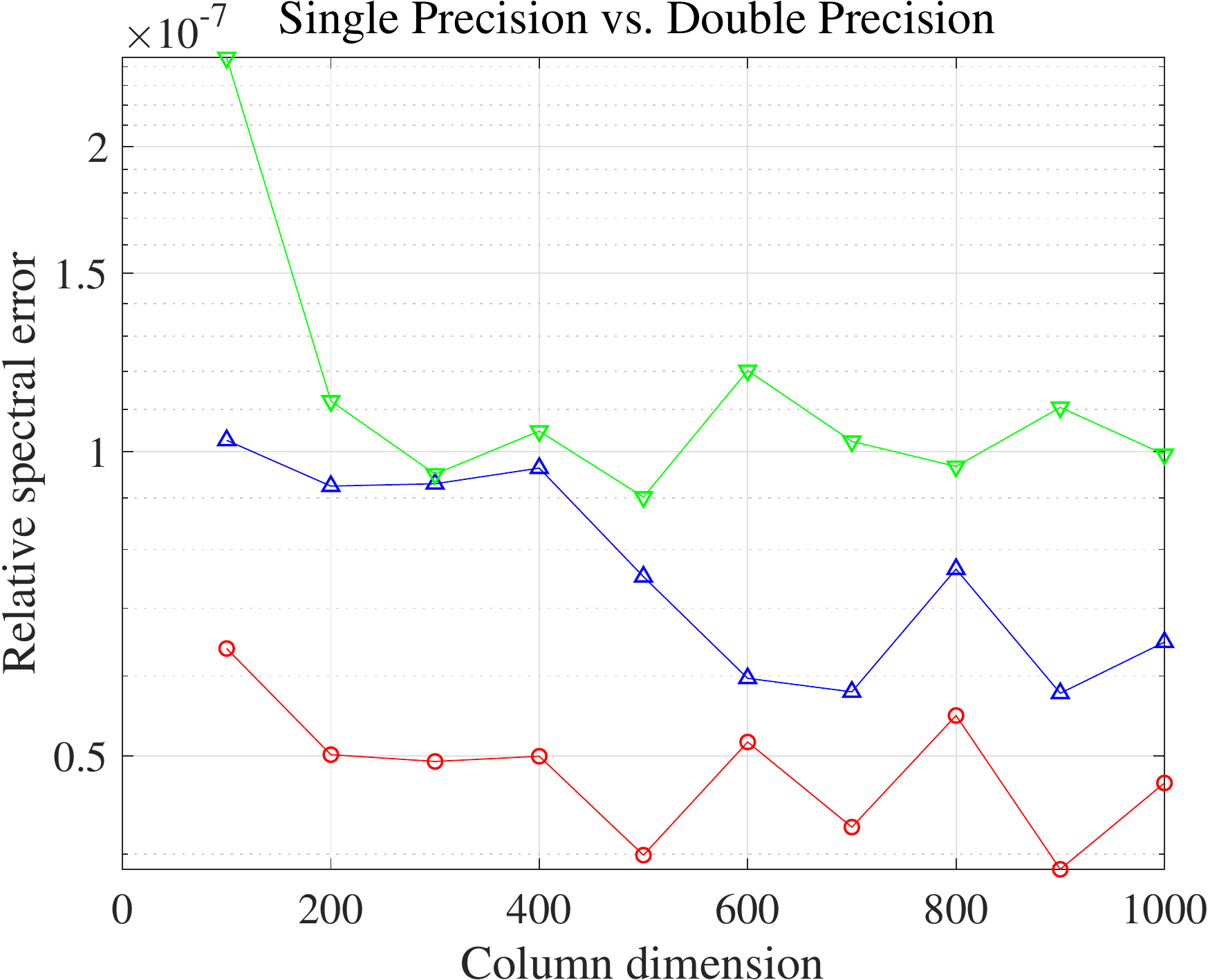}
\caption{Relative spectral error of single precision ID with respect to target rank (left) and column dimension (right).}
\label{fig:lp_dp}
\end{figure}
When we execute mixed single precision ID, we are still using the columns from the double precision matrix in order to construct our column skeleton. On the other hand, in single precision ID, the columns come from the single precision matrix. Generally, we should expect this to introduce error into our approximation. However, our three test matrices have value ranges well within single precision round-off and therefore our single precision ID should still perform well in these tests. This is confirmed by the left panel of Figure~\ref{fig:lp_dp}, as our scheme matches the results generated by the double precision algorithm to within unit round-off for almost all target rank values for the medium and fast decay datasets. Again, the slow decay data matrix is more of a challenge for our scheme, though the relative errors, all of which are less than $10^{-6}$, are acceptable given that the low-rank approximation itself introduces far more significant error in all three cases. 

Our single precision ID is robust with respect to increases in column dimension, as shown in the right panel of Figure~\ref{fig:lp_dp}. In these test cases, single precision ID is almost exactly as close to the double precision scheme as mixed single precision ID. Let us consider the case in which we approximate all $1000$ columns of the matrix. In the slow decay test case, our optimal matrix $2$-norm error for a rank-$20$ approximation is $\sigma_{21} = 1/21$~\cite{eckart1936approximation}, so the discrepancies between the two methods are negligible compared to the low-rank approximation error. In the medium decay test case, the optimal error is given by $\sigma_{21} = 1/(21)^2 = 1/441$, and again the discrepancy is negligible. Finally, in the fast decay test case, the optimal approximation error is roughly $5 \times 10^{-6}$, indicating that single precision ID accumulates round-off error which is dominated by the low-rank approximation error in all three datasets.

%

We now conduct the same set of numerical experiments, but in simulated half precision, meaning that our data matrices are cast to half precision, but round-off error is accumulated in single precision. Because unit round-off in half precision is significantly larger than in single precision ($4.9 \times 10^{-4}$ as opposed to $6 \times 10^{-8}$), we expect much worse performance in some cases, including the possibility of our algorithms breaking as a result of numerical underflow.
\begin{figure}
\centering
\includegraphics[width=0.49\linewidth]{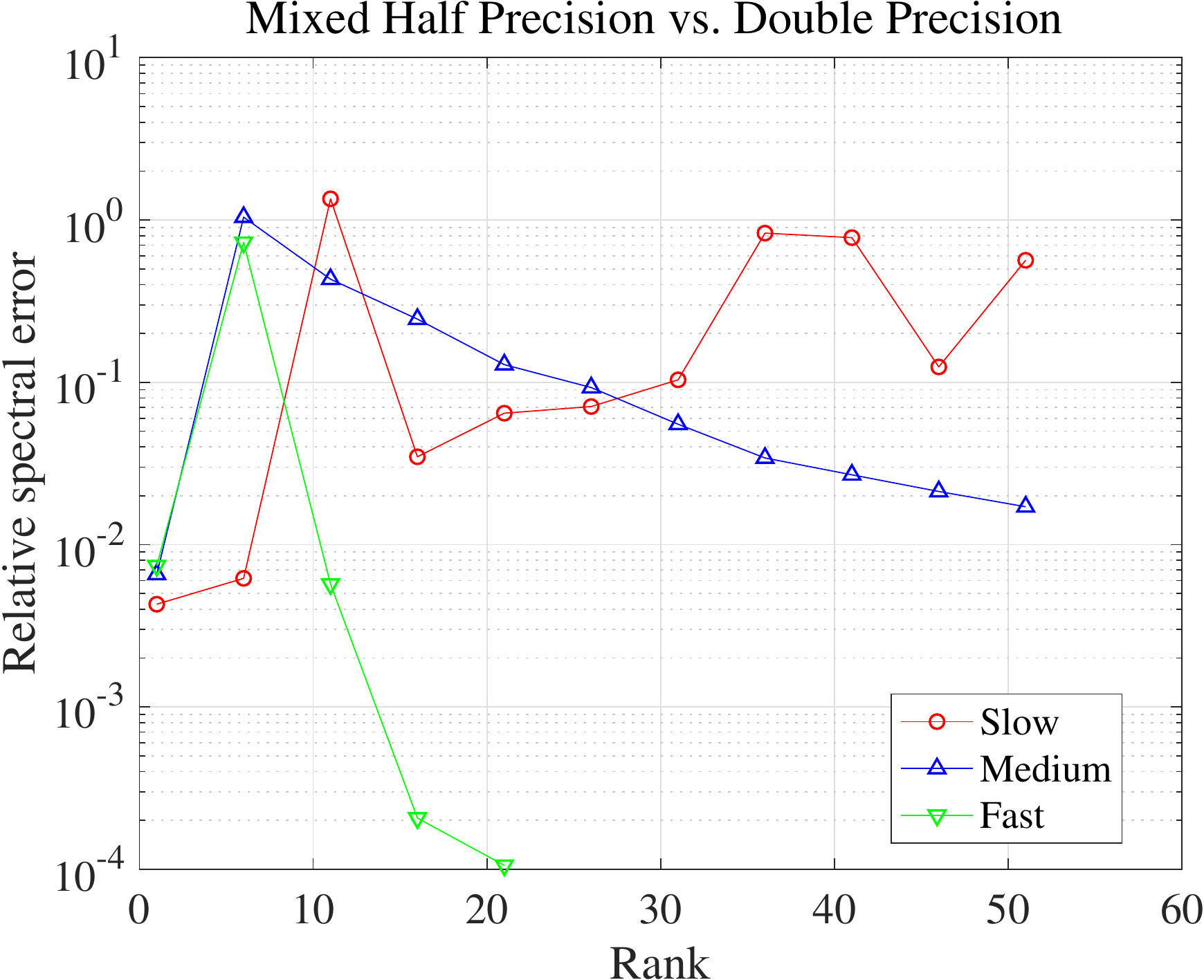}
\includegraphics[width=0.49\linewidth]{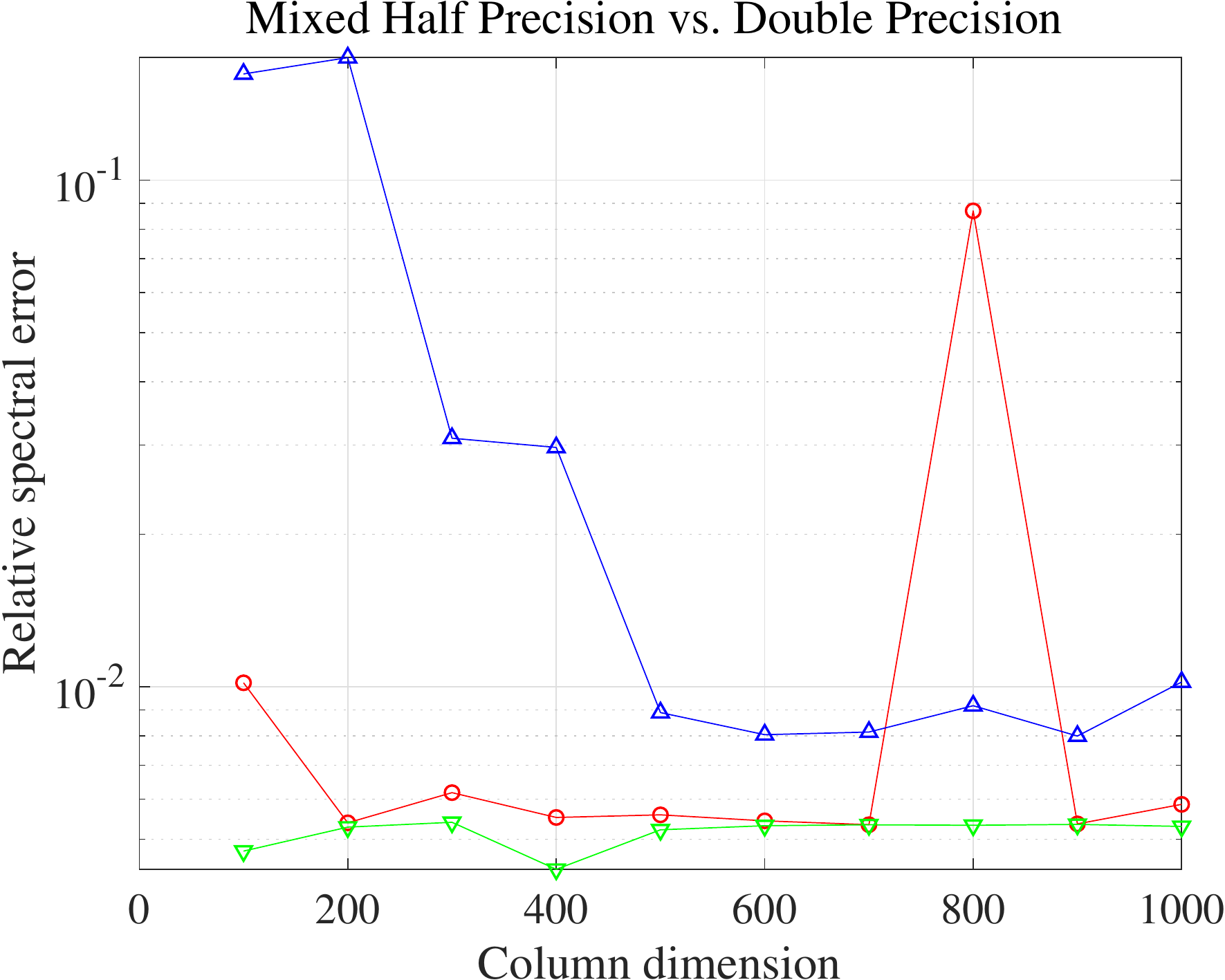}
\caption{Relative spectral error of mixed half precision ID with respect to target rank (left) and column dimension 
(right). The algorithm breaks on the fast decay matrix for target ranks exceeding 21 due to numerical underflow.}
\label{fig:mp_dp_hp}
\end{figure}

We observe a significant difference in performance between mixed half precision ID and double precision ID. In the left panel of Figure~\ref{fig:mp_dp_hp}, we see that our scheme performs reasonably well on the medium decay and fast decay test cases for low target rank values. For higher target rank values, the method breaks in the fast decay case due to numerical underflow, while struggling on the slow decay test matrix. In the right panel of Figure~\ref{fig:mp_dp_hp} we observe that we are within 1-2 digits of accuracy of the double precision approximation in almost all cases, with the algorithm performing the best on the fast decay test matrix and worst in the medium decay case. Mixed half precision ID is able to generate reasonably accurate low-rank approximations for all three matrices, and is therefore useful in situations where computational bottlenecks may necessitate the use of half precision arithmetic. 
\begin{figure}
\centering
\includegraphics[width=0.49\linewidth]{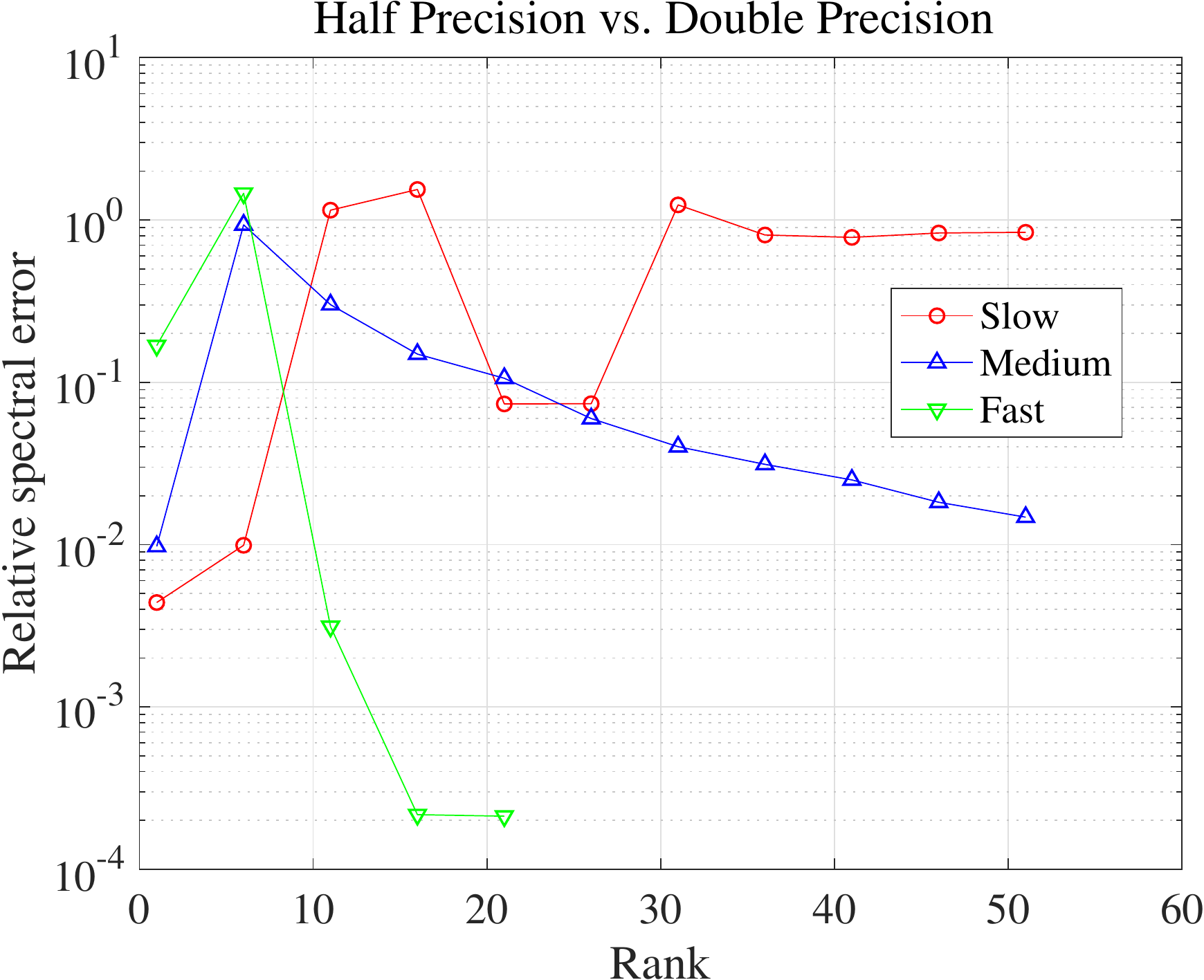}
\includegraphics[width=0.49\linewidth]{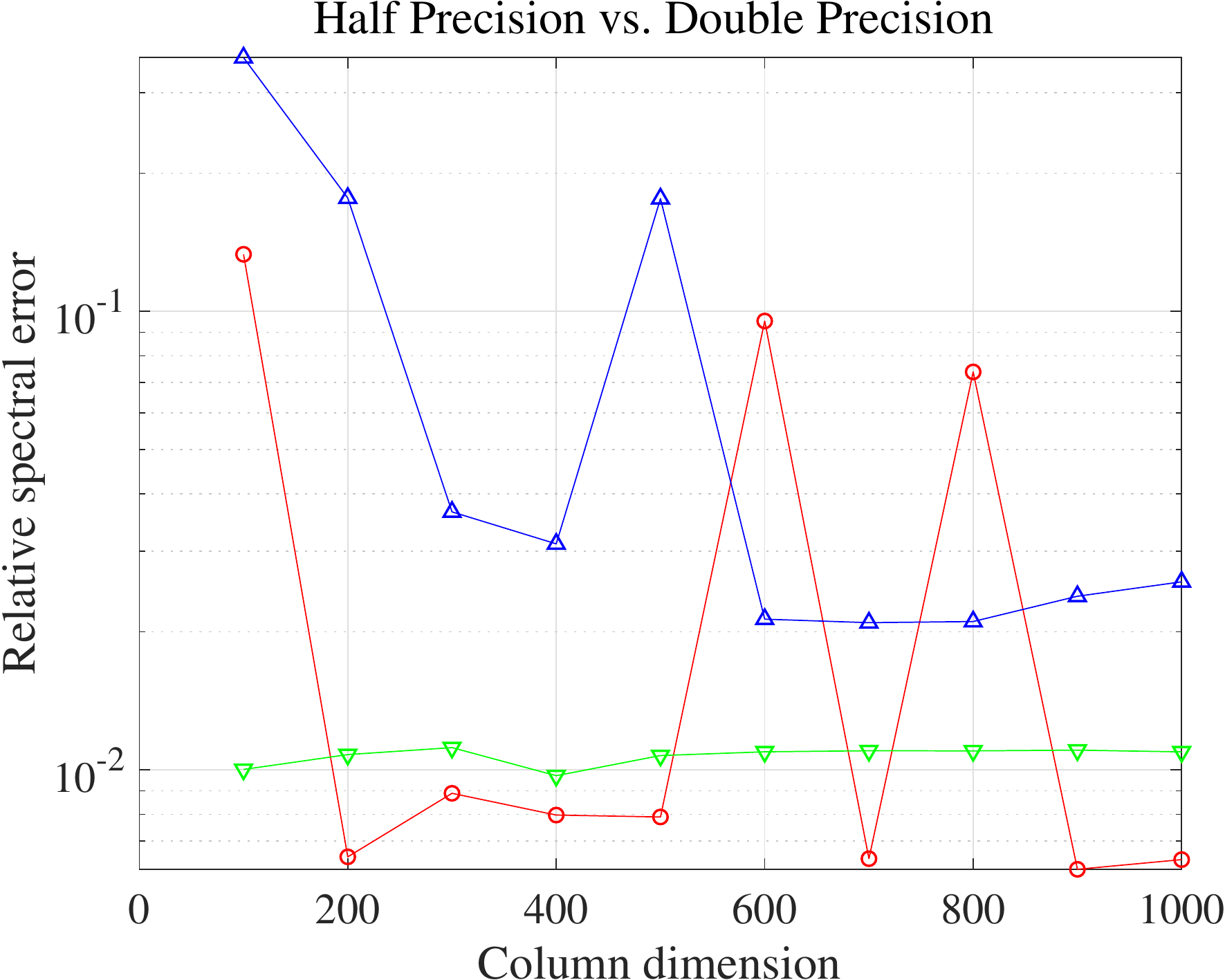}
\caption{Relative spectral error of half precision ID with respect to target rank (left) and column dimension 
(right). The algorithm breaks on the fast decay matrix for target ranks exceeding 21 due to numerical underflow.}
\label{fig:lp_dp_hp}
\end{figure}

 Finally, we compare half precision ID to double precision ID. In this test case, we expect results similar to those between mixed half precision and double precision, which is verified by Figure~\ref{fig:lp_dp_hp}. In the left panel of the figure, we see that our method struggles in the slow decay test cases, improves as rank is increased in the medium decay test case, and performs quite well in the fast decay test case for low target rank values before breaking, again due to numerical underflow. As we increase the column dimension of our matrix, the half precision algorithm maintains about 1 digit of accuracy for all three test matrices, performing best in the fast decay test case and far less consistently in the slow and medium decay test cases (right panel of Figure~\ref{fig:lp_dp_hp}).

\subsection{Comparison to ground truth}

\begin{figure}
\centering
\includegraphics[width=0.49\linewidth]{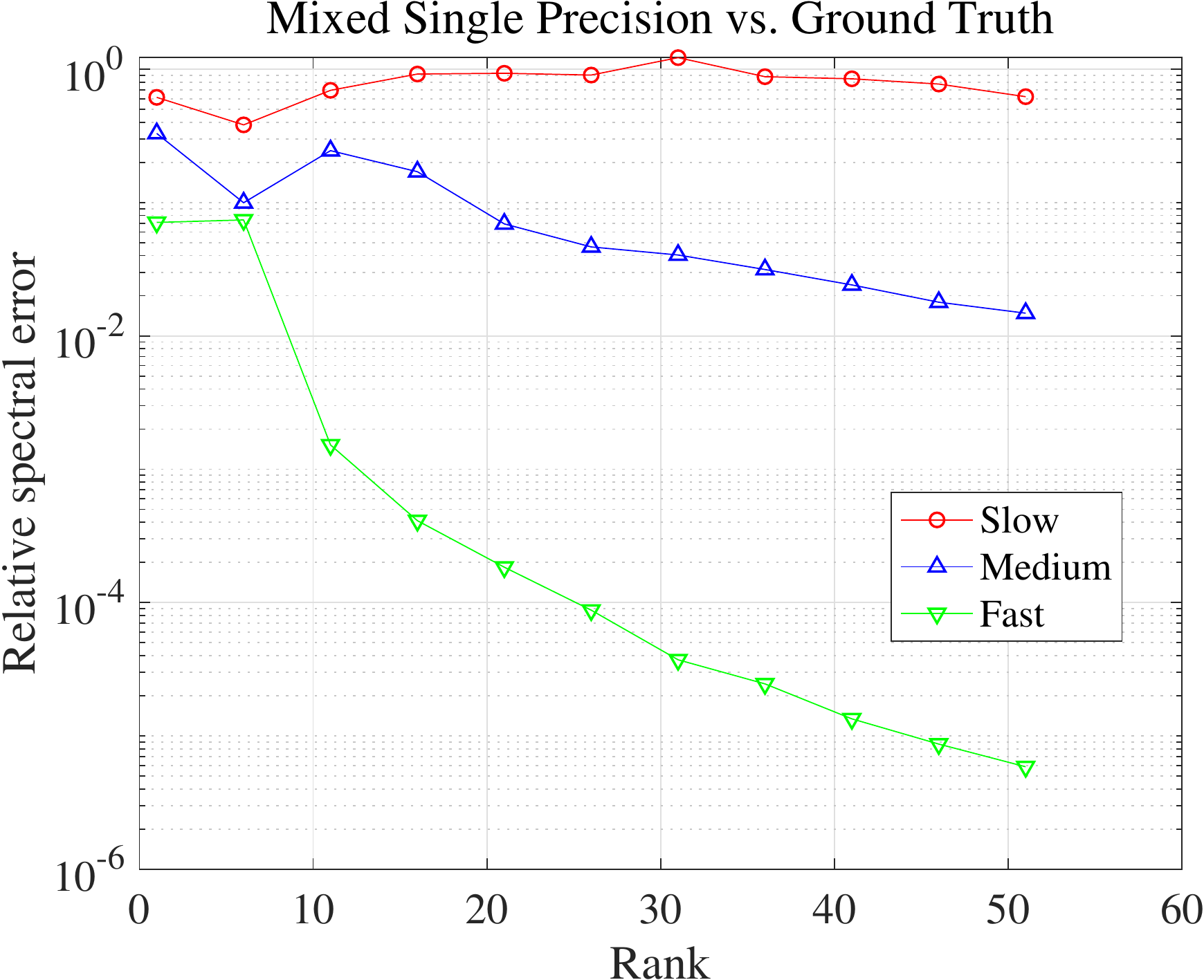}
\includegraphics[width=0.49\linewidth]{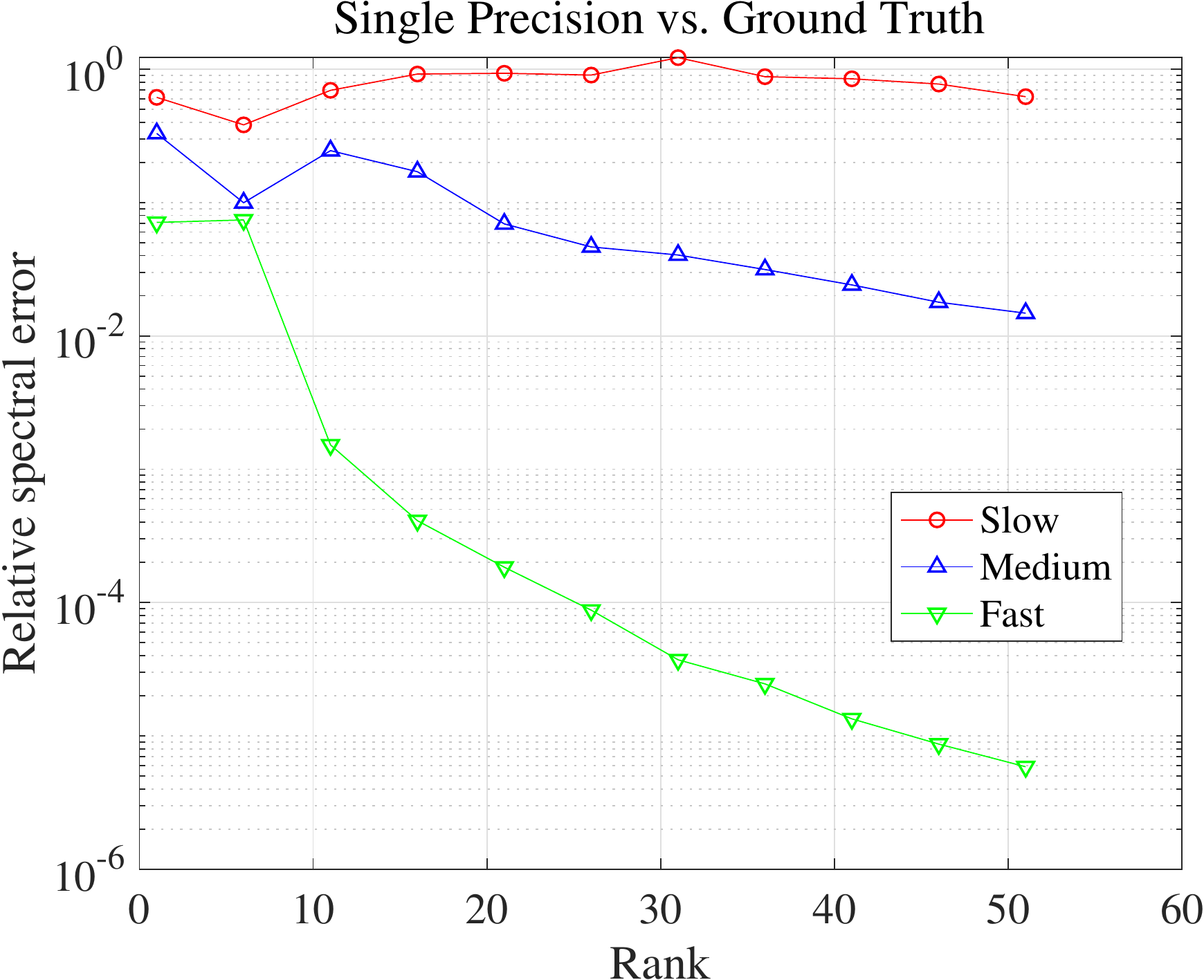}
\caption{Relative spectral error of mixed single precision ID (left) and single precision ID (right) compared to the ground truth double precision matrix with respect to target rank (left) and column dimension (right).}
\label{fig:mp_truth}
\end{figure}
We now measure the performance of our single precision algorithms against the ground truth double precision matrix. In Figure~\ref{fig:mp_truth}, we observe that for all three test datasets, our algorithm achieves improvements in error as target rank is increased.  We also evaluate the approximation error of single precision ID relative to the ground truth matrix. Because a single precision ID approximation is generated using a single precision column subset as well as a single precision coefficient matrix, we should expect some loss of accuracy. Comparing the two panels of Figure~\ref{fig:mp_truth}, we see that single precision ID performs almost exactly as well as mixed precision ID, demonstrating that at least for the three test matrices we have selected for our experiments, single precision ID is an effective tool for low rank matrix approximation.

\begin{figure}
\centering
\includegraphics[width=0.49\linewidth]{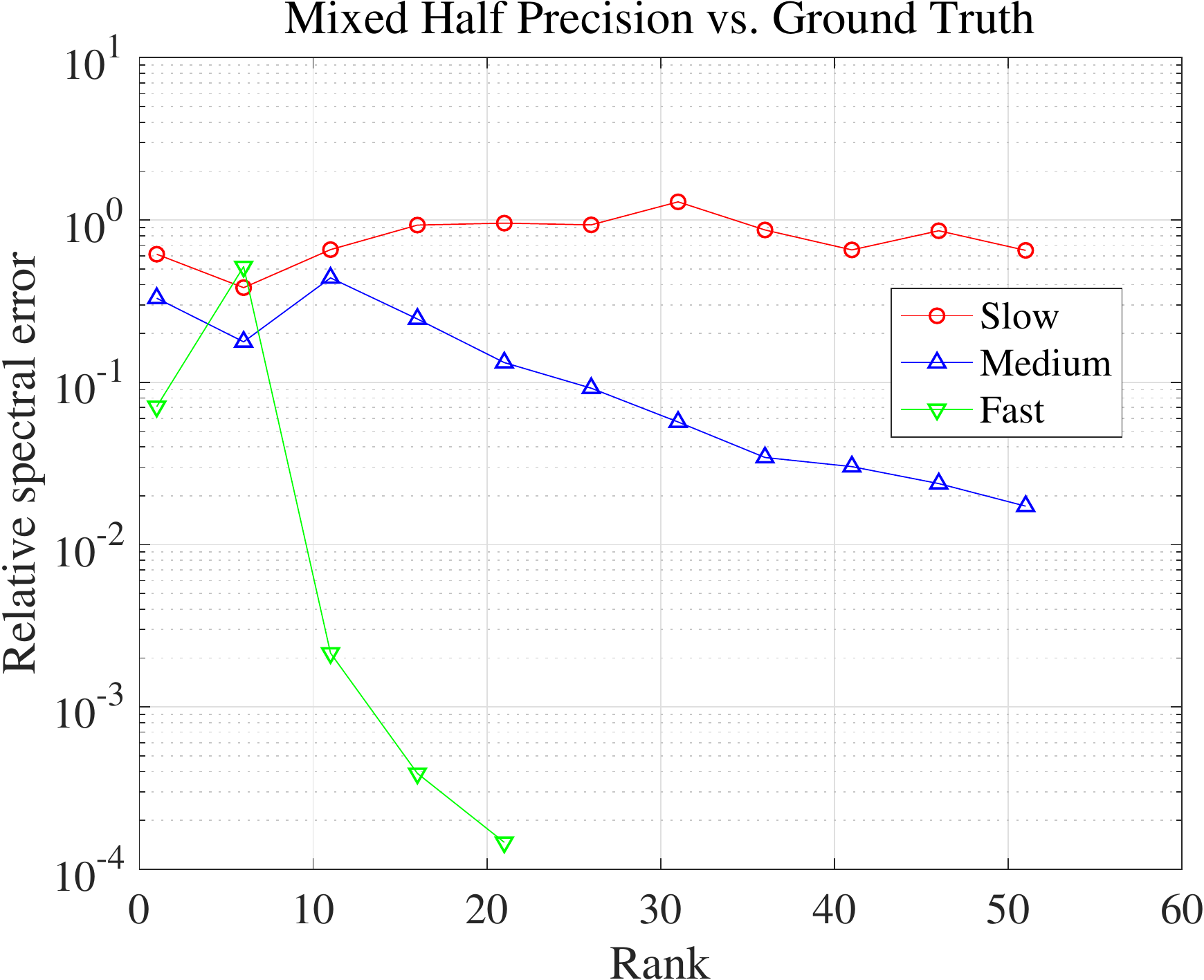}
\includegraphics[width=0.49\linewidth]{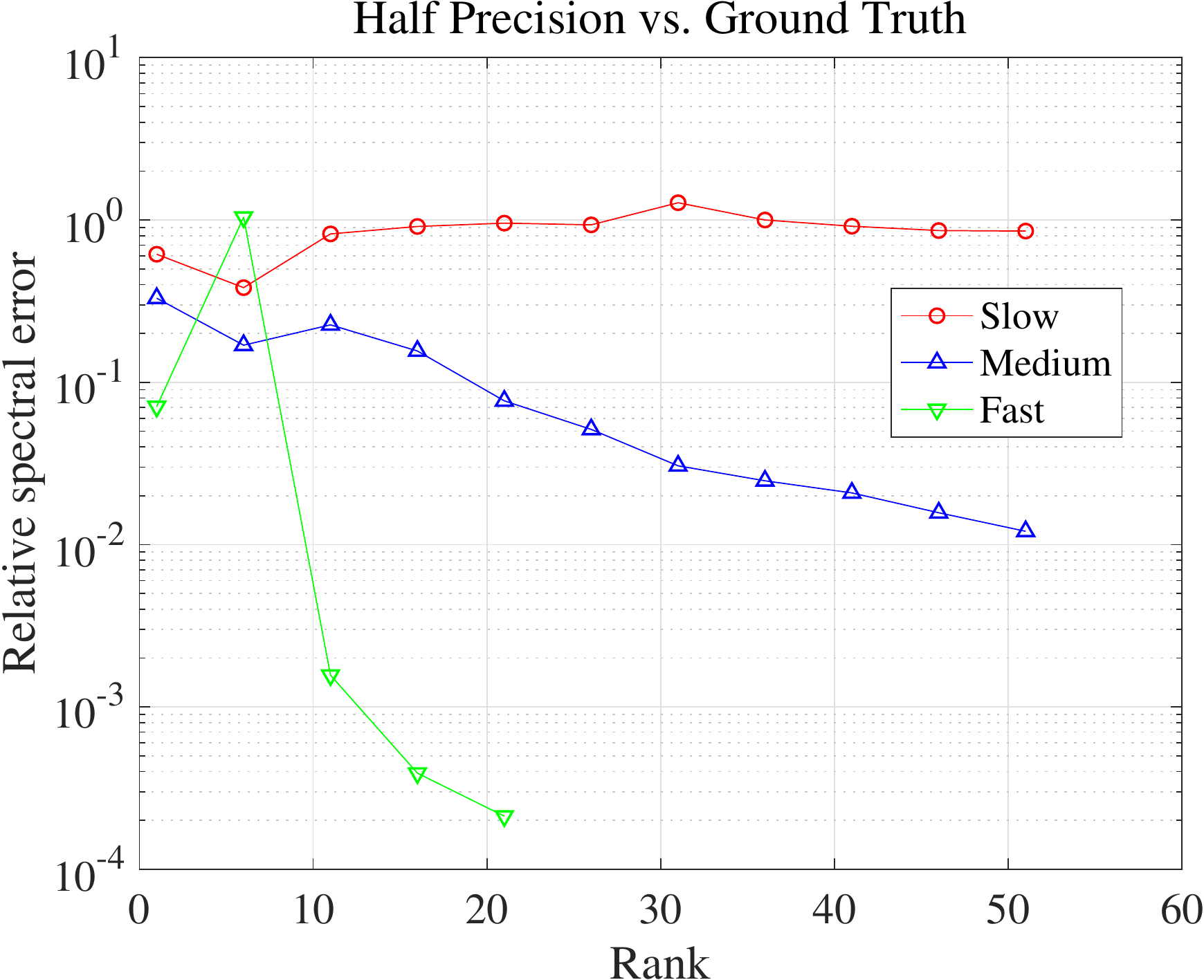}
\caption{Relative spectral error of mixed precision ID in mixed half (left panel) and half precision (right panel) compared to the ground truth double precision matrix with respect to target rank.}
\label{fig:mp_truth_hp}
\end{figure}
We now conduct the same set of numerical experiments, but in simulated half precision. We observe the failure of our scheme to approximate the fast decay data matrix for any target rank exceeding $21$ in Figure~\ref{fig:mp_truth_hp}. 
We also see that mixed half precision ID performs about as well as single precision ID, though mixed half precision ID breaks on the fast decay data matrix for target ranks exceeding $21$. We compare our entirely low half precision ID approximation to the ground truth double precision data matrix in the right panel. In this test, we expect significant errors due to the fact that we are approximating our original double precision matrix with an interpolative decomposition comprised of half precision columns, as well as a half precision coefficient matrix. We see that our algorithm struggles in the slow decay test case, performs reasonably well in the medium decay case, and gives good approximations to our fast decay matrix before breaking for too large of target rank values.

\begin{figure}[t]
\centering
  \includegraphics[width=\linewidth]{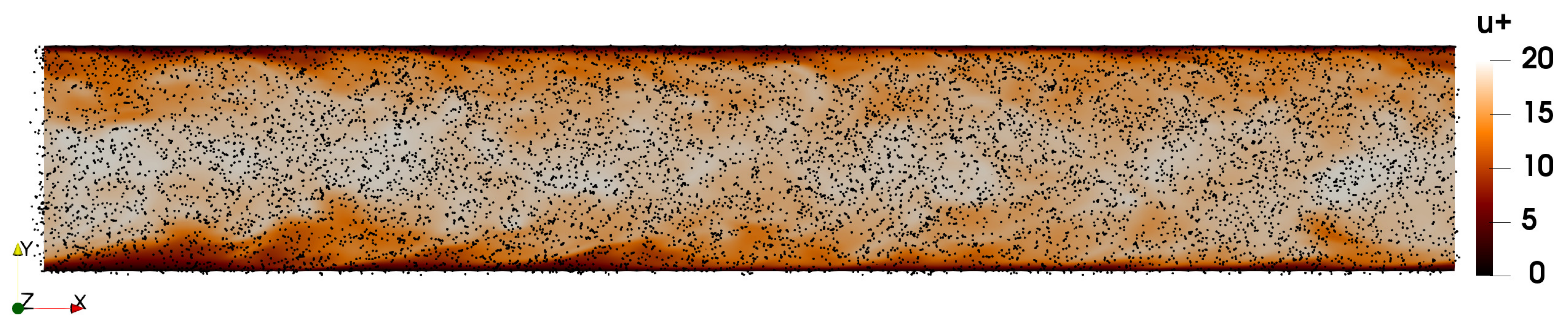}
  \caption{Cross section of a particle-laden turbulent channel flow~\cite{dunton2019pass}.}
  \label{fig:turb_channel_flow}
\end{figure}

\section{Reduced order modeling for particle-laden turbulence}
\label{sec:particles}
As a real-world dataset, we select stream-wise velocity data taken from particles suspended in a turbulent channel flow (see Figure~\ref{fig:turb_channel_flow}) with Stokes number $St^+ = 1$ (see~\cite{jofre2018multi,dunton2019pass} for background on the problem setup and corresponding dataset). We seek to use our mixed and low precision algorithms to (1) identify the `most important' subset of particles in a large scale simulation and (2) use the information from these particles to compute statistics of the system of interest. In this application, our dataset will be the stream-wise velocities of 5000 particles measured over 10000 time steps in a turbulent channel flow with periodic boundary conditions. We will use the mixed precision ID algorithms to select the most important subset of these particles, and then use this subset of particles to predict the time evolution of the entire system.

\begin{table}[h]
\centering
\begin{tabular}{|p{2.3cm}|p{2.3cm}|p{2.3cm}|p{2.3cm}| }
\hline
Dimensions & $\sigma_{50}/\sigma_{1}$ & $\sigma_{n}/\sigma_{1}$ & value range\\
\hline
$10000 \times 5000$ & $5.0 \times 10^{-4}$ & $3.3\times 10^{-9}$ & $4.1 \times 10^{6}$\\
\hline
\end{tabular}
\caption{Characteristics of particle stream-wise velocity data matrix.}
\label{tab:particledataproperties}
\end{table}

\begin{figure}
\centering
\includegraphics[width=0.49\linewidth]{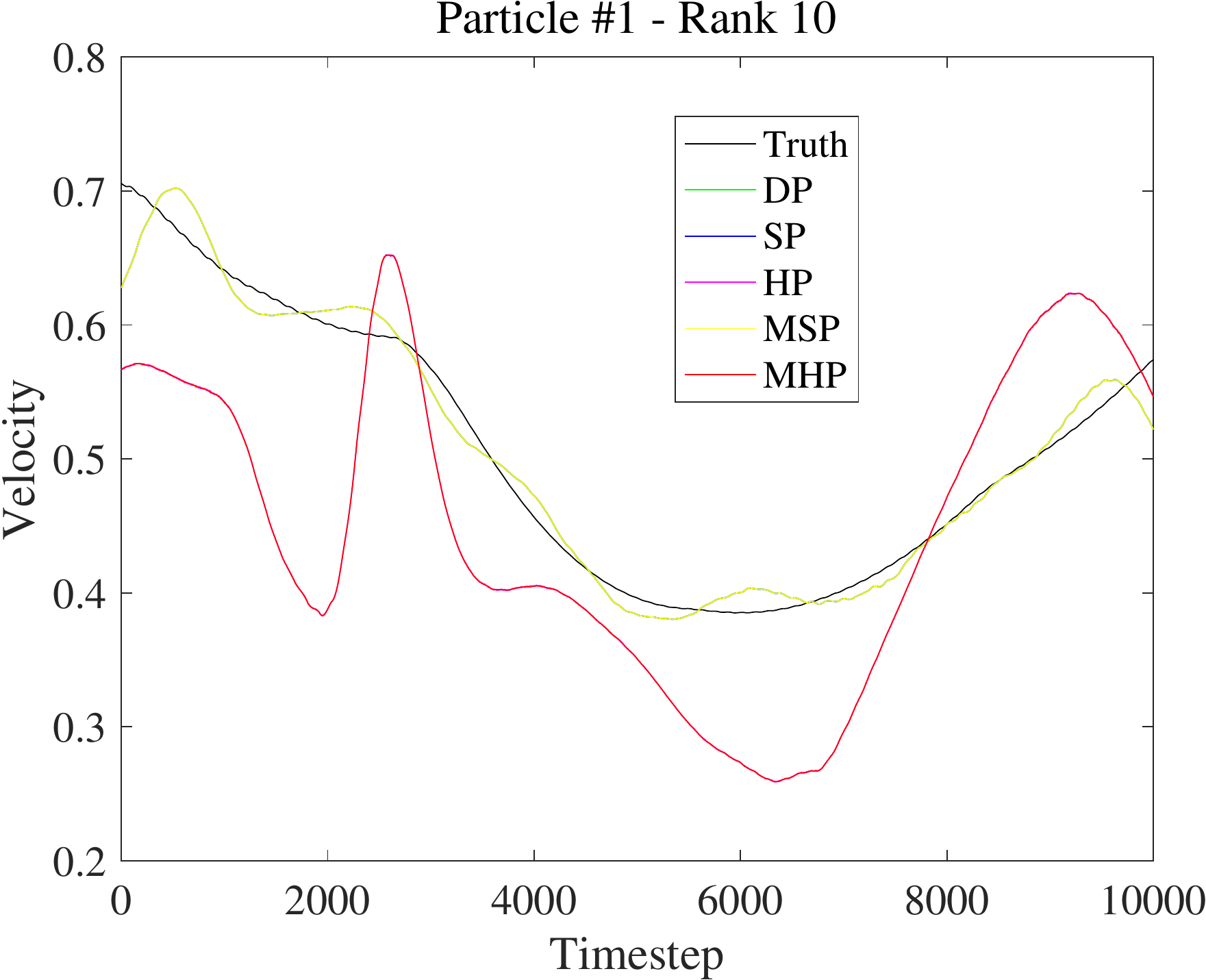}
\includegraphics[width=0.49\linewidth]{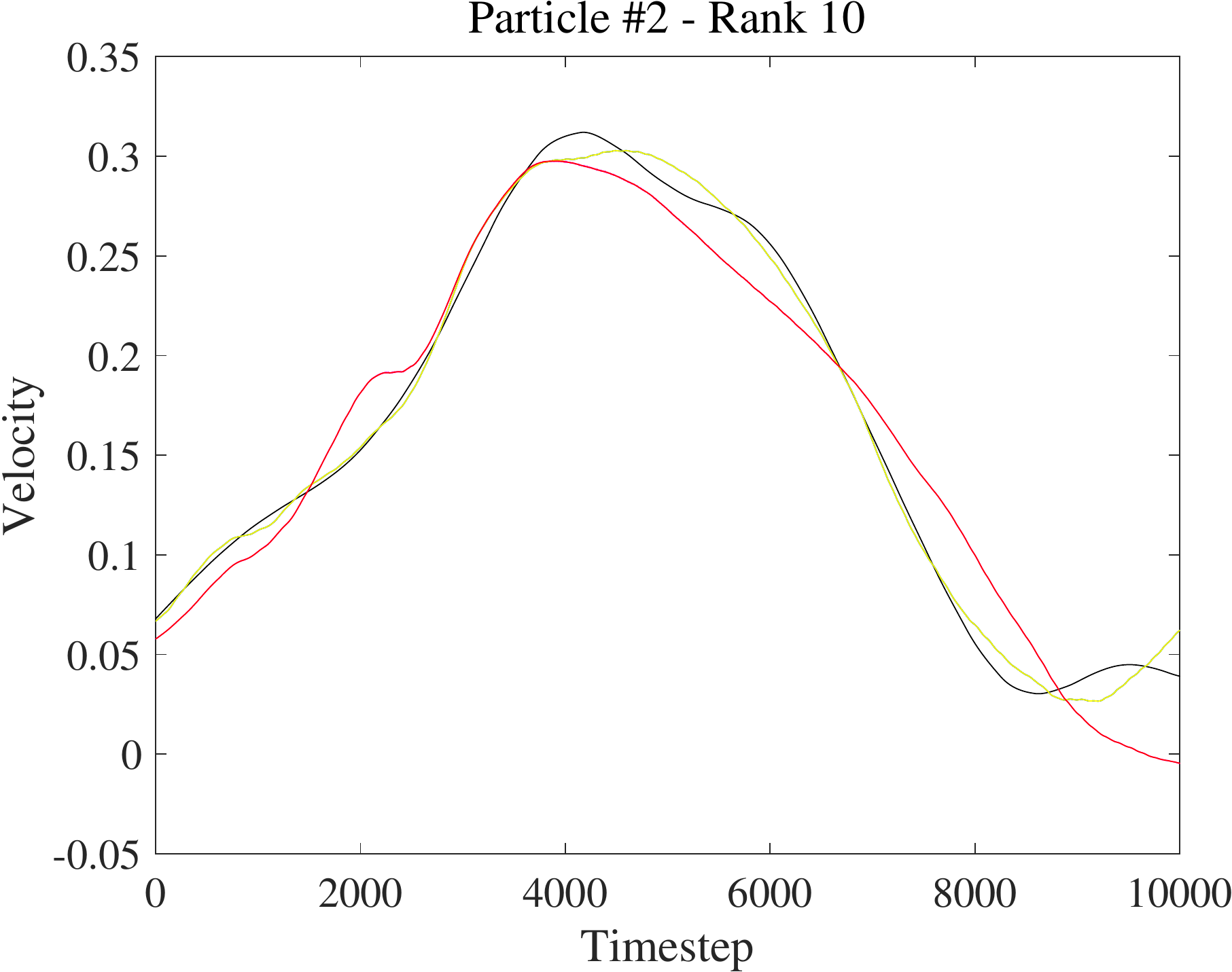}
\includegraphics[width=0.49\linewidth]{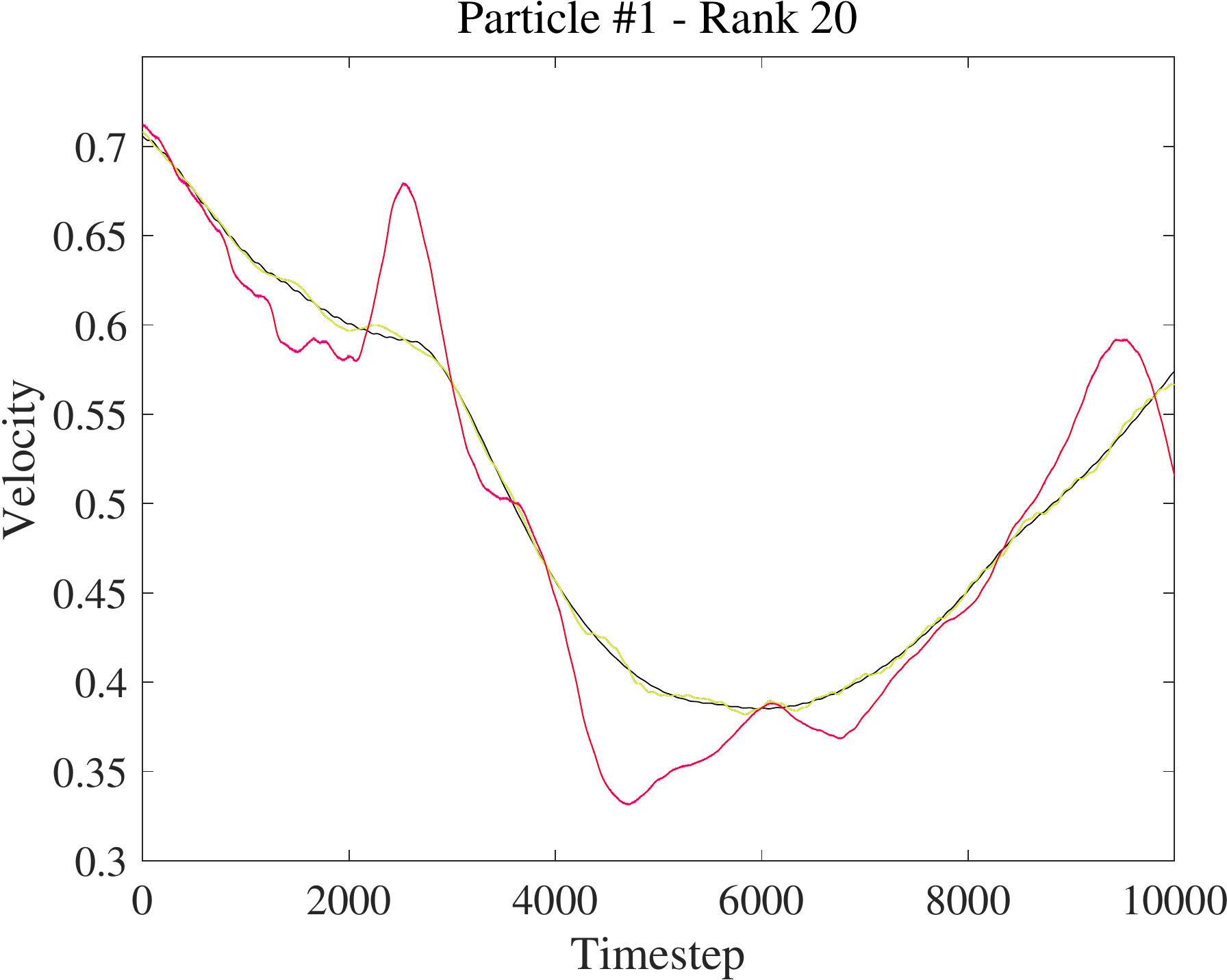}
\includegraphics[width=0.49\linewidth]{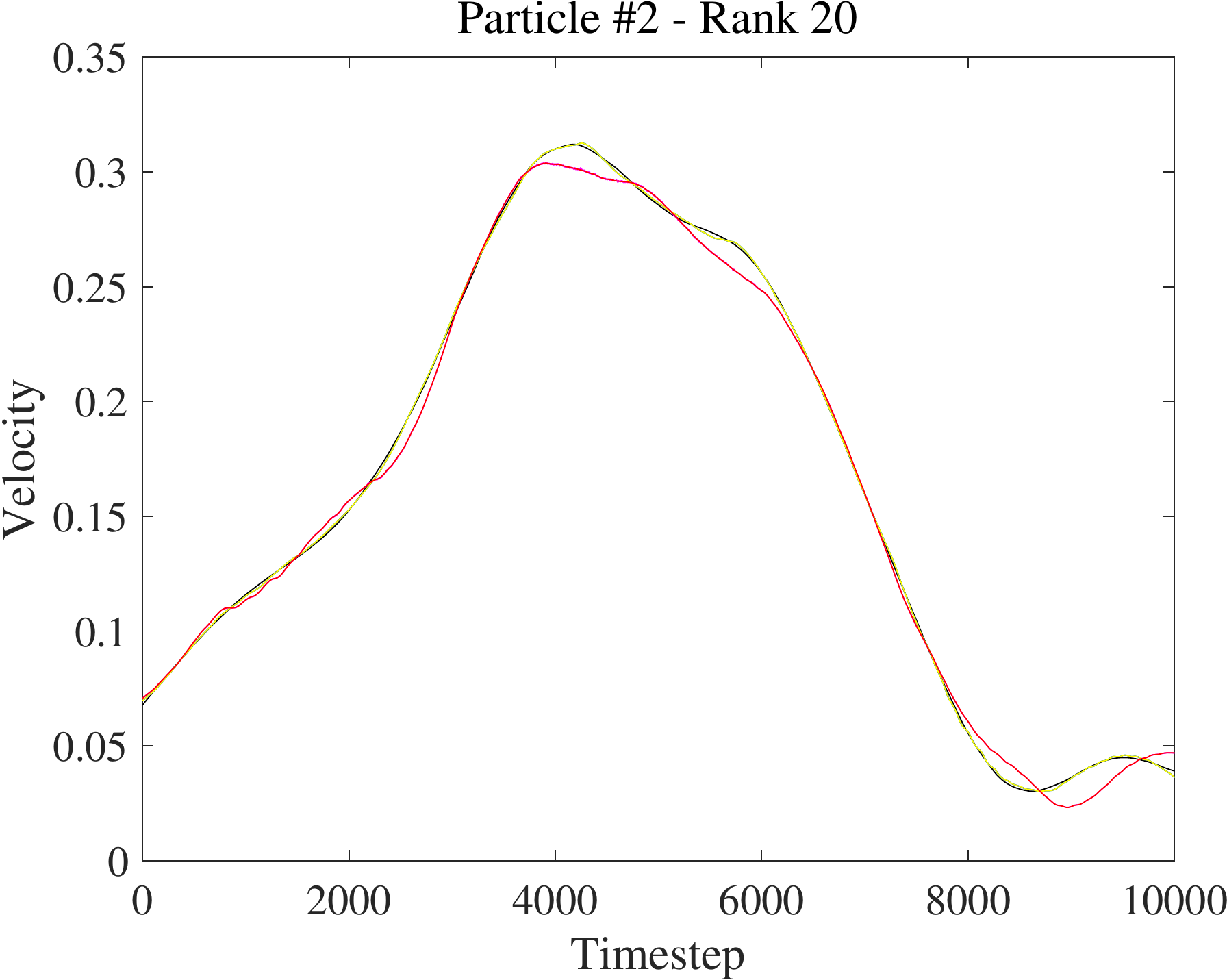}
\includegraphics[width=0.49\linewidth]{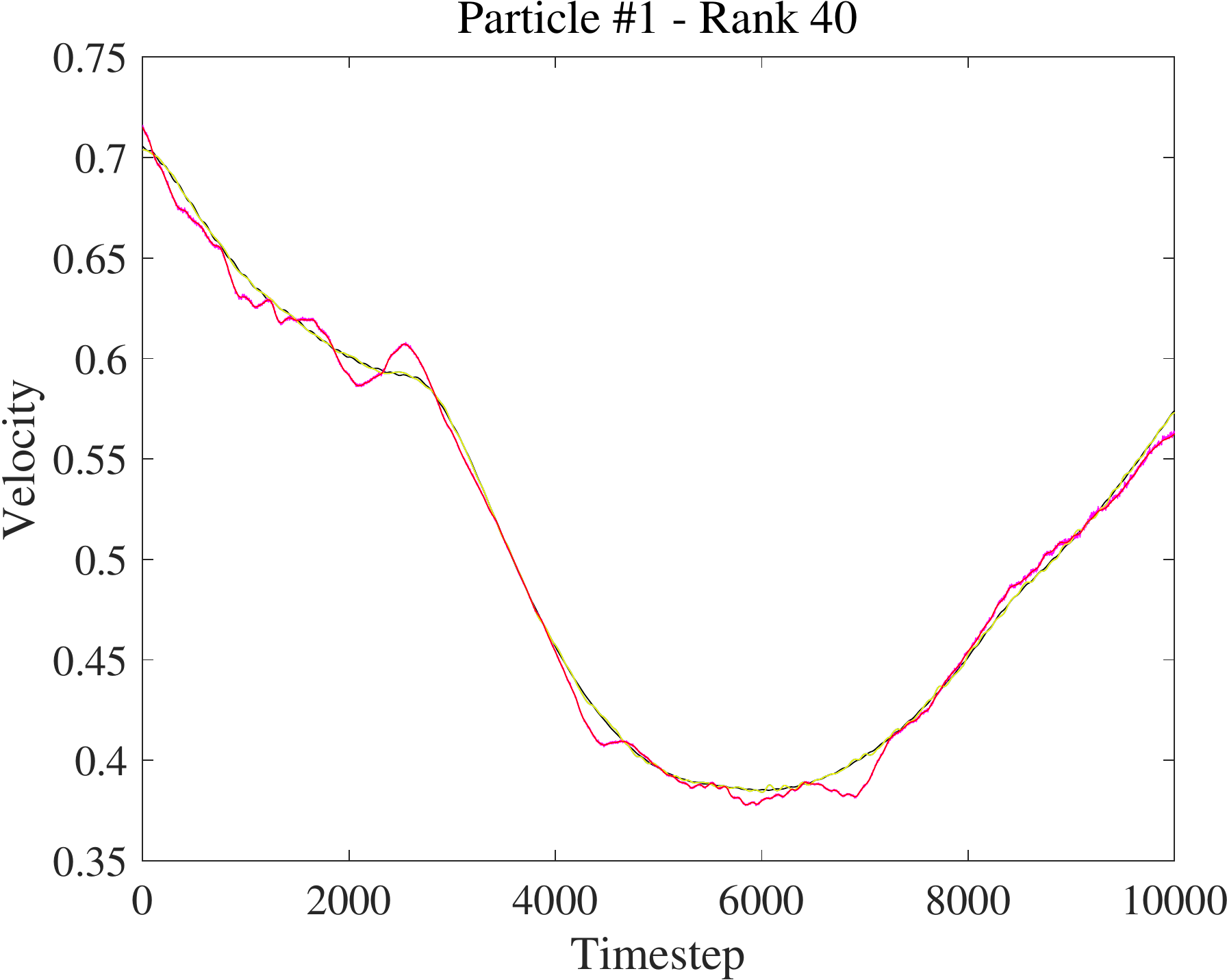}
\includegraphics[width=0.49\linewidth]{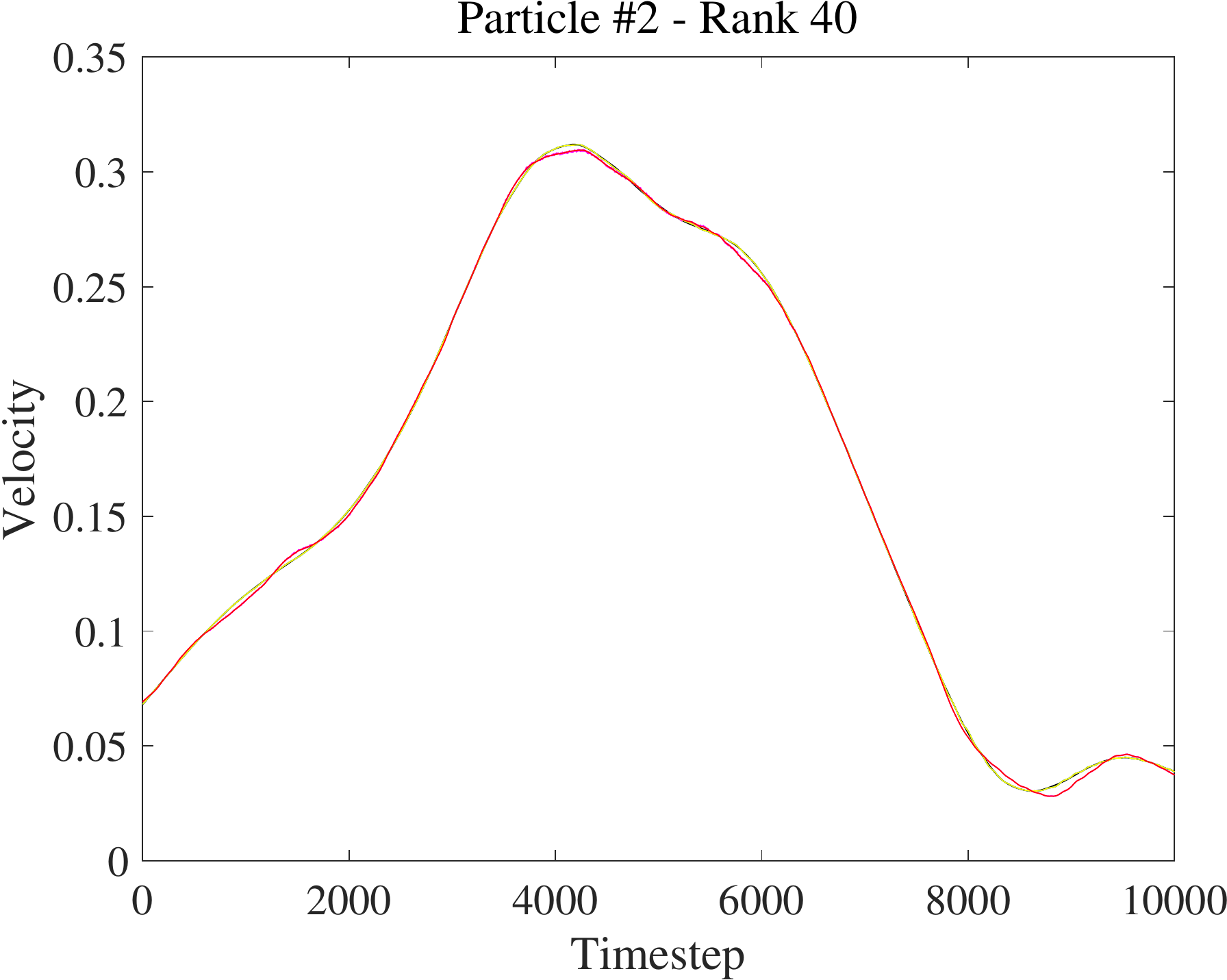}
\caption{Velocity evolution over $10000$ time-steps for two particles approximated using double precision, single precision, half precision, mixed single precision, and mixed half precision ID for target rank values of $10,20,40$.}
\label{fig:particle_laden_velocities}
\end{figure}

\begin{table}
\centering
\begin{tabular}{ |p{3.5cm}||p{2.5cm}|p{2.5cm}|p{2.5cm}| p{2.5cm}| p{2.5cm}|  }
 \hline
 \multicolumn{4}{|c|}{Particle 1 (\# 149) velocity evolution mean-squared approximation error} \\
 \hline
 ID Algorithm & Target Rank 10 & Target Rank 20 & Target Rank 40\\
 \hline
Double Precision  & 2.2e-2   & 5.3e-4 &  6.8e-5\\
Single Precision & 2.2e-2 & 5.3e-4   & 6.8e-5\\
Half Precision & 9.8e-1 & 1.1e-1 & 3.9e-3\\
Mixed Single Precision  & 2.2e-2 & 5.3e-4 & 6.8e-5\\
Mixed Half Precision& 9.8e-1 & 1.1e-1 & 3.9e-3\\
 \hline
\end{tabular}
\caption{Mean-square approximation error of particle data using our five low and mixed precision ID algorithm for three target rank values of $10,20$, and $40$. Errors are computed for the particles which are not selected by the algorithm for the reduced order model (in which case the approximation is exact by definition of the ID).}
\label{tab:part1}
\end{table}

\begin{table}

\centering
\begin{tabular}{ |p{3.5cm}||p{2.5cm}|p{2.5cm}|p{2.5cm}| p{2.5cm}| p{2.5cm}|  }
 \hline
 \multicolumn{4}{|c|}{Particle 2 (\# 2064) velocity evolution mean-squared approximation error} \\
 \hline
 ID Algorithm & Target Rank 10 & Target Rank 20 & Target Rank 40\\
 \hline
Double Precision  & 4.6e-3  & 7.2e-5 &  6.4e-6 \\
Single Precision & 4.6e-3 & 7.2e-5   & 6.4e-6\\
Half Precision & 5.2e-2 & 3.2e-3 & 2.7e-4 \\
Mixed Single Precision  & 4.6e-3 & 7.2e-5 & 6.4e-6\\
Mixed Half Precision& 5.2e-2 & 3.2e-3 & 2.7e-4\\
 \hline
\end{tabular}
\caption{Mean-square approximation error of particle data using our five low and mixed precision ID algorithm for three target rank values of $10,20$, and $40$. Errors are computed for the particles which are not selected by the algorithm for the reduced order model (in which case the approximation is exact by definition of the ID).}
\label{tab:part2}
\end{table}
\label{sec:partladenturbchannelflow}
We use all four variations of mixed and low precision ID, as well as a regular double precision ID, in order to construct our reduced order models. We use subsets of $10$, $20$, and $40$ particles to predict the velocities of two particles (not originally chosen by our column selection) over the $10000$ time steps for which they are tracked in our flow. The ground truth trajectories and the five approximate trajectories generated using double precision, single precision, half precision, mixed single precision, and mixed half precision ID are shown in Figure~\ref{fig:particle_laden_velocities}. We expect that as we increase our target rank, our predicted trajectories should become more accurate. This is supported by our results for both particles - as we increase the target rank (number of particles in our reduced order model) we see that in both cases all five algorithms converge to the ground truth. We also expect that the double and single precision methods will perform superior to half precision, which is reflected in the plots - the worst performing methods of the five are the half and mixed half precision IDs. 

In order to further distinguish the performance of the five methods, we compute the mean squared approximation error over the $10000$ time-steps for which we trace our two particles. We expect that the double precision method will perform the best of the five, while both single precision algorithms should perform almost as well. Finally, our half precision schemes are expected to incur the most error. This is confirmed by the mean squared errors for the two particles reported in Tables~\ref{tab:part1} and~\ref{tab:part2}. The single precision algorithms both perform as well as the double precision for both particles and all three target rank values. The half precision method perform the worst and as poorly as one another. This demonstrates that even if we use the double precision columns from the ground truth matrix to construct our approximation, in some cases, the round-off error in the coefficient matrix $\bm{P}$ dominates the error when we compute an ID in half precision.

\section{Conclusions}
\label{sec:conc}
In this work we derive mixed and low precision algorithms for computing the matrix interpolative decomposition. Experiments demonstrate that our schemes are effective in single precision. Low-rank approximation is not as well suited to half precision arithmetic, therefore, our methods do not perform as well in half precision, but nevertheless succeed in generating low-rank approximations in some test cases. Half precision appears to be most effective for matrices with sharp singular value decay and small target rank values.  We apply our algorithms to a reduced order modeling problem, using them to model entire systems of particles suspended in turbulent flow using only a fraction of the entire system. An immediate next step of this work is to derive representative error estimates using, e.g., probabilistic error analysis~\cite{higham2018new}, as well as to incorporate other ideas from numerical linear algebra such as matrix sketching into our algorithms for increased speed. Moreover, blocked implementations of our algorithms will help mitigate roundoff accumulation, particularly in half precision, which will accommodate the latest developments in half precision hardware.
\bibliographystyle{siamplain}
\bibliography{Ref}

\end{document}